\documentclass[conference,letterpaper]{IEEEtran}
\usepackage[letterpaper, left=0.75in, right=0.75in, bottom=0.75in, top=0.75in]{geometry}
\def\IEEEtitletopspace{18pt}
\IEEEoverridecommandlockouts
\usepackage{array}
\usepackage{hyperref}
\usepackage{cite}
\usepackage{amsmath,amssymb,amsfonts}
\usepackage{algorithmic}
\usepackage{comment}
\usepackage{graphicx}
\usepackage{adjustbox}

\usepackage{textcomp}
\usepackage{breqn}
\usepackage{amsmath}
\usepackage{xcolor}
\usepackage{booktabs}
\usepackage{xcolor}
\usepackage{subcaption}
\usepackage{caption}
\newcount\Comments  
\usepackage{color}
\definecolor{darkgreen}{rgb}{0,0.5,0}
\definecolor{purple}{rgb}{1,0,1}
\definecolor{teal}{rgb}{0,0.4627,0.5804}
\newcommand{\kibitz}[2]{\ifnum\Comments=1\textcolor{#1}{#2}\fi}

\def\BibTeX{{\rm B\kern-.05em{\sc i\kern-.025em b}\kern-.08em
    T\kern-.1667em\lower.7ex\hbox{E}\kern-.125emX}}
\begin{document}

\title{Flight Path Optimization with Optimal Control Method\\

%


\thanks{
$^{1}$G. Su and K. Liu are with the University of California at Berkeley, Berkeley, CA 94702, USA. 

$^{2}$X. Cheng is with the University of Illinois at Chicago, Chicago, IL 60607, USA.

$^{3}$S. Feng (corresponding author) is with the Hong Kong University of Science and Technology, Clear Water Bay, Hong Kong.

$^{4}$J. Song is with the University of Toronto, Toronto, Ontario M5S 1A1 Canada.

$^{5}$J. Chen is with the University of British Columbia, Vancouver, BC, Canada V6T 1Z2.

$^{6}$C. Zhu is with the Tsinghua University, Beijing, 100084, P. R. China.

$^{7}$H. Lin is with Northwestern University, Evanston, IL, 60208, USA.

}
}

\author{Gaofeng Su$^{1}$, Xi Cheng$^{2}$, Siyuan Feng$^{3,*}$, Ke Liu$^{1}$, Jilin Song$^{4}$, Jianan Chen$^{5}$, Chen Zhu$^{6}$, Hui Lin$^{7}$}

\maketitle

\begin{abstract}

This paper is based on a crucial issue in the aviation world: how to optimize the trajectory and controls given to the
aircraft in order to optimize flight time and fuel consumption. This study aims to provide elements of a response to
this problem and to define, under certain simplifying assumptions, an optimal response, using Constrained Finite
Time Optimal Control (CFTOC). The first step is to define the dynamic model of the aircraft in accordance with the
controllable inputs and wind disturbances. Then we will identify a precise objective in terms of optimization and
implement an optimization program to solve it under the circumstances of simulated real flight situation. Finally,
the optimization result is validated and discussed by different scenarios.

\end{abstract}

\begin{IEEEkeywords}
Constrained Finite Time Optimal Control; Flight Path Optimization; Aircraft Point Mass Model
\end{IEEEkeywords}

\section{Introduction}
In today’s world, we rely on various modes of transportation, from ground \cite{cheng2024autonomous},\cite{shafiee2024using,liu2023physics} to air \cite{liu2023airborne,liu2021miles}, to travel from one place to another. These transportation methods contribute significantly to CO2 emissions, exacerbating global warming and the consequent rise in sea levels \cite{cheng2024carsharing, cheng2024electric}. Because of
this, the transportation industry has seen a major shift in the technology used just to lower these emissions \cite{cheng2023estimating,cheng2024using}. But this
has not been the same case across all different modes of transportation, especially in the airplane industry. Airplanes
contribute to 12:5\% of the global CO2 pollution and 80\% of the emissions are from flights of over 1500km. This is one of the main reasons this project focuses on flight path optimization. Flight paths contribute to the amount of
fuel used, distance travelled which taking into account the wind speed, air pressure and distance between the two
locations. They can provide a lot of advantages such as:
1. By optimizing the flight path to consume the least fuel while finding the shortest path possible, we can
decrease emissions by a large factor. This reduces the overall environmental impact of airplanes.
2. A decrease in amount of fuel used saves the airlines a lot of money while benefiting the passengers by
decreasing the cost of their overall journey.
3. They can minimize the impact of the weather on the aircraft thus reducing the amount of maintenance
required and saving the airlines more money.
4. Better management of Air Traffic due to predetermined routes.

There are a lot of challenges that come with creating flight plans which include:
1. Taking into account the weather over the whole journey. Weather here refers to pressure of air and wind
speed at all points during the journey.
2. The weight and amount of fuel being carried. This can also be related to the weight and size of the airplane
as larger amounts of fuel are usually carried by larger planes
There are many ways to design and determine flight paths. The method that we chose is Constrained Finite Time
Optimal Control. A major reason for this was that there was not a lot of research about using Constrained Finite
Time Optimal Control (CFTOC) to design flight paths. One reason might be that CFTOC is mostly used in Model
Predicted Control (MPC) with receding horizon approach. Hence, it is interesting to discuss the possibility of using
CFTOC for a nonlinear, large scale optimization problem, which we have further explored through this paper.

This paper is organized as follows. Section II reviews related work in aircraft trajectory optimization. Section III describes our multi-criteria optimization framework, while Section IV details its implementation. Section V presents case study results, showcasing the framework's impact on efficiency and sustainability. Section VI discusses broader implications for air traffic management and environmental policy. Section VII concludes with a summary and future research directions, highlighting the potential for advancements in sustainable aviation practices.

\section{Literature Review}
This section explores recent inquiries into flight path optimization, specifically considering weather variations and wind disturbances. These studies aim to refine aircraft trajectory optimization under diverse atmospheric conditions by merging technological advancements with environmental foresights.

The groundwork laid by Fechner and Schmehl (2018) \cite{fechner2018flight} introduces an innovative optimization technique for the flight trajectories of tethered wings in airborne wind energy systems. Their methodology integrates advanced wind models and reeling strategies to enhance efficiency and power output under turbulent conditions, enabling more effective wind energy exploitation in aviation. Building upon this, Lindner et al. (2020) \cite{lindner2020flight} present a pioneering real-time optimization framework for aircraft flight paths that utilizes weather forecast-based Corridors of Optimization (CoO). Their approach not only achieves significant fuel savings but also improves flight efficiency, striking a balance between operational adaptability and air traffic predictability.

Franco et al. (2010) \cite{franco2010minimum} investigate fuel consumption optimization at constant altitudes by introducing a variable-Mach cruise strategy and highlighting the impact of arrival timings on fuel efficiency. Girardet et al. (2013, 2014) \cite{girardet2013generating} \cite{girardet2014wind} propose methodologies for optimizing aircraft routes to minimize travel time and congestion. Their contributions leverage favorable winds and employ sophisticated algorithms to sidestep aerial hazards, marking significant advancements in dynamic and weather-responsive flight planning.

Ng et al. (2014) \cite{ng2014optimizing} develop an innovative algorithm that incorporates wind conditions into aircraft trajectory optimization, resulting in notable fuel savings from Anchorage, Alaska. Meanwhile, Botkin et al. (2017) \cite{botkin2017aircraft} explore aircraft control strategies under windshear using differential game theory and viability theory, enhancing flight safety under adverse weather conditions. Selecký et al. (2013) \cite{selecky2013wind} refine UAV trajectory planning in windy scenarios with an improved Accelerated A* algorithm, improving flight path accuracy and collision avoidance.

Rosenow et al. (2021) \cite{rosenow2021advanced} assess the impact of in-flight trajectory optimization on fuel efficiency and airspace capacity, advocating for the seamless integration of Trajectory-Based Operations (TBOs) into current flight operations. Valenzuela and Rivas (2014) \cite{valenzuela2014optimization} examine air traffic control constraints on cruise flight optimization, emphasizing the need for flexible optimization strategies within existing regulatory frameworks.

Another paper goes beyond existing research by presenting a multi-criteria optimization framework that addresses fuel efficiency, flight duration, and environmental impacts simultaneously in aircraft trajectory planning \cite{majumder2016flight}. Unlike singularly focused studies, our unified model amalgamates these considerations. Authors also explore the application of real-time data and machine learning predictions for dynamic flight path adjustments \cite{yu2019flight}, ushering in a new era for static optimization models that currently dominate the field.

Deori et al. (2015) \cite{deori2015model} develop a Model Predictive Control model for aircraft path following, resonating with our optimization challenge. Their exploration of appropriate aircraft models and constraints sheds light on our methodology. However, our investigation differs in that we utilize Constrained Finite Time Optimal Control to identify optimal flight paths, with a specific focus on travel time or fuel cost. This synthesis not only engages with existing literature but also propels the field towards innovative methodologies for enhancing air travel efficiency and sustainability.

\section{Study Objective}
Implement Constrained Finite Time Optimal Control (CFTOC) to optimized the flight path with minimum travel 
time focus/minimum fuel cost focus from Chicago O’Hare International Airport (ORD) to San Francisco International
Airport (SFO) taking forecast wind data into account.

\section{Technical Description}
The implementation and resolution of the problem is based on four parts. First of all, the dynamic model of the
aircraft is configured and developed, according to the objectives set and the degree of precision required. In a
second step, we do the same for the wind, and integrate this modeling with the aircraft one. Once the dynamic
model is established, it must then be discretized, which can make it easier and more feasible to be programmed
and solved. Indeed, given the complexity and non-linearity of the function describing the dynamic behavior of
the system, an analytical solution to our problem is diffcult to find. And so we choose to solve the problem
numerically. Finally, the final optimization problem will be implemented in the last phase.

\subsection{Aircraft model}
A aircraft point mass model (PMM) will be introduced in the path planning optimization problem. The original
model is in 3-dimension, however, to simplify the problem, the original model was reduced into a 2-dimension
model.

\begin{equation}
\begin{aligned}
\dot{x} &= v \cos({\theta}) + w_{x} \\
\dot{y} &= v \sin({\theta}) + w_{y}  \\
\dot{v} &= \frac{T-D}{m}=\frac{2 T - C_{d} \rho A v^{2}}{2m} \\
\dot{m} &= - \eta T \\
\dot{\theta}&= \frac{g}{v}tan(\phi)\\
\end{aligned}
\end{equation}

\vspace{1\baselineskip}
\noindent
Where,\\
$x$: west-east distance [m]\\
$y$: north-south distance [m]\\
$w_{x}$: wind disturbance on x direction [m/s]\\
$w_{y}$: wind disturbance on y direction [m/s]\\
$v$: aircraft ground speed [m/s]\\
$m$: aircraft mass [kg]\\
$\theta$: heading angle [rad]\\
$C_{d}$: drag coefficient [-]\\
$\rho$: air density [kg/$m^3$]\\
$A$: aircraft wing area [$m^2$]\\
$\eta$: thrust specific fuel consumption coefficient [kg/N$\cdot$s]\\
$T$: aircraft thrust force [N]\\
$\phi$: bank angle [rad]\\

In this model, $\dot{\theta}$ involves the aircraft speed, which appears in the denominator. This implies that v cannot be 0.
For a path following problem, this is acceptable. However, in a path planning problem, nonzero speed restriction
implies that the aircraft can only reach the destination at the end of the given time horizon (i.e. the aircraft cannot
reach the destination with a minimum travel time and stop there). Also, for nonlinear programming, the nonlinear
solver is very sensitive to the form of the optimization problem and the constraints parameters. Hence, it is better
to avoid putting speed in the denominator.
Therefore, the expression of $\dot{\theta}$ was reformulated and replaced by a new input variable $\phi$ [rad/s], which represent
the turning rate. The final aircraft point mass model that will be used in path planing:

\begin{equation}
\begin{aligned}
\dot{x} &= v \cos({\theta}) + w_{x} \\
\dot{y} &= v \sin({\theta}) + w_{y}  \\
\dot{v} &= \frac{2 T - C_{d} \rho A v^{2}}{2m} \\
\dot{m} &= - \eta T \\
\dot{\theta}&= \varphi\\
\end{aligned}
\end{equation}

\subsection{Wind model}
The Historical Track database contains records of 1336 flights between Chicago O'Hare International Airport (ORD) and San Francisco International Airport (SFO) from July 2013 to August 2013. Information such as departure time, geographic coordinates during the flight, and aircraft identification (ACID) can be retrieved from this database. These flights were operated using two types of aircraft: the Airbus A319 and the Airbus A320.

Upon analyzing the paths of these 1336 flights, we identified five commonly traversed Air Route Traffic Control Centers (ARTCC Zones) along the route from ORD to SFO: ZAU, ZDV, ZLC, ZMP, and ZOA. We will use the boundaries of these zones to segment wind data. An example of past flight tracks, like VRD211, is depicted in the figure below.

In our modeling approach (MPC model), we rely on forecasted wind data sourced from the National Center for Environmental Prediction (NCEP). NCEP generates these predictions using the Rapid Refresh (RAP) method, an hourly-updated system managed by the National Oceanic and Atmospheric Administration (NOAA). The RAP system comprises a numerical forecast model and an analysis/assimilation system. The wind data provided are in two-dimensional coordinates, including longitude and latitude for location and wind speed components for both x and y directions.

The wind data is stored in grib2 format, requiring a grib file reader such as pygrib for data retrieval in our project.

An issue with the RAP forecast data is its utilization of a hybrid vertical coordinate system rather than an isobaric coordinate system. Converting between these coordinate systems is complex. To address this, we calculated the average pressure for each layer in the forecast model and selected the layer with the closest pressure to 250mb, corresponding to the typical cruise altitude of selected flights (approximately 32,000 feet).

\begin{figure}
\begin{center}
    \includegraphics[scale=0.58]{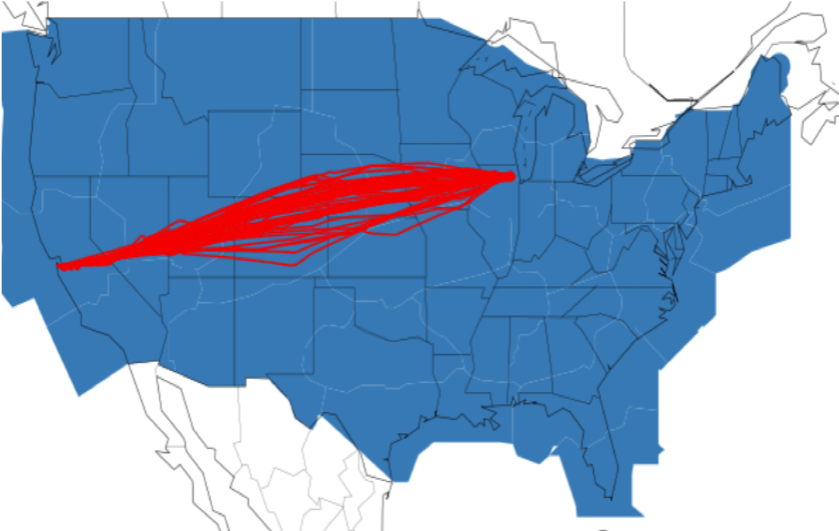}
    \caption{Actual Flight Tracks from ORD to SFO}
    \end{center}
\end{figure}

\begin{figure}
\begin{center}
    \includegraphics[scale=0.21]{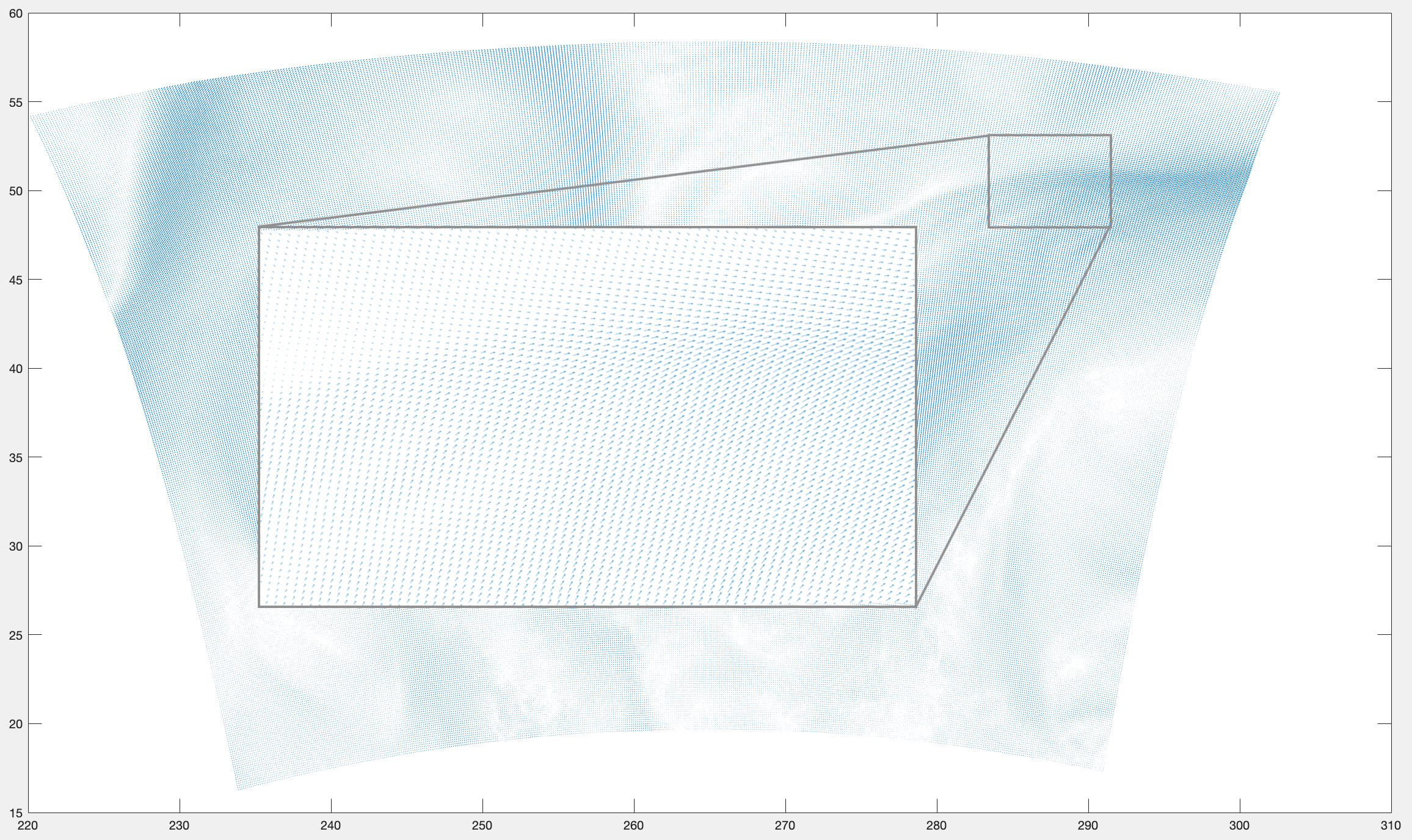}
    \caption{Vector field generated with the raw wind data collected on 07/02/13 at 00:00.}
    \end{center}
\end{figure}

Wind components are added into our dynamics to create a more realistic MPC model.

\begin{equation}
\begin{aligned}
\dot{x_k} &= v_k \cos({\theta_k}) + w_{x_k} \\
\dot{y_k} &= v_k \sin({\theta_k}) + w_{y_k}  \\
\end{aligned}
\end{equation}

However, the raw wind data obtained from NCEP cannot be directly combined with our dynamics since the raw wind data components are not directly compatible with the velocity of the aircraft. To convert the raw wind data into the compatible form with the velocity dynamics of the aircraft, we developed a mathematical expression.

We first pick out the five regions defined by Federal Aviation Administration that all the real aircraft path goes through. The wind force are considered in two perpendicular directions u and v. Then we take the average of the wind force for all the 5 time slots that is accessible, and obtained the average wind force spatial distribution in the two directions u and v. We also transform the 2D coordinate systems from WGS84 to WGS 84/Pseudo-Mercator with python package pyproj, and then apply another linear transformation to set the origin of the coordinates system to the departure location. The unit is set to kilometer. 

However, the discrete wind data will be hard to handled as input into the aircraft dynamics. Thus, we construct two polynomial functions to fit the two wind force spatial distribution.Here, we use x and y to represent coordinates along longitude and latitude directions. The trick is how to determine the degree and items of the polynomial functions. From the plots above, we could find the wind force in u direction is high correlated in x and y, while the correlation is weak for the v direction. Based on this finding, we try different polynomial functions with degrees ranging from 3 to 9 to see their performance in fitting the real wind force distribution. A least square method is applied to specify the parameter of the polynomial functions and the best performances are found when the degree is set to 5. The results are shown below. Figure 2 to Figure 4 show decent fitting performance of the two designed polynomial functions.

\vspace{1\baselineskip}

$ w_{x_k} = a_1 y_k^4 + a_2 y_k^3 + a_3 y_k^2 + a_4 y_k + a_5 + a_6 x_k^4 + a_7 x_k^3 + a_8 x_k^2 + a_9 x_k + a_{10} + a_{11} y_k^3 x_k + a_{12} y_k^2 x_k^2 + a_{13} y_k x_k^3 $

\vspace{1\baselineskip}

where
\[
\begin{pmatrix}
\ a_1 \\
\ a_2 \\
\ a_3 \\
\ a_4\\
\ a_5 \\ 
\ a_6 \\ 
\ a_7 \\ 
\ a_8 \\ 
\ a_9 \\ 
\ a_{10} \\ 
\ a_{11} \\ 
\ a_{12} \\ 
\ a_{13} \\ 
\end{pmatrix}
\approx
\begin{pmatrix}
\ 5.404 \cdot 10^{-12} \\
\ -7.525 \cdot 10^{-9} \\
\ -1.010 \cdot 10^{-5} \\
\ 1.8023 \cdot 10^{-3}\\
\ 3.054 \cdot 10^{-1}\\ 
\ 1.071 \cdot 10^{-12} \\ 
\ 8.131 \cdot 10^{-9} \\ 
\ 1.957 \cdot 10^{-5} \\ 
\ 1.360 \cdot 10^{-2} \\ 
\ 3.054 \cdot 10^{-1} \\ 
\ -4.493 \cdot 10^{-13} \\ 
\ 1.372 \cdot 10^{-12} \\ 
\ -1.971 \cdot 10^{-12} \\ 
\end{pmatrix}
\]

\vspace{1\baselineskip}

$ w_{y_k} = b_1 y_k^4 + b_2 y_k^3 + b_3 y_k^2 + b_4 y_k + b_5 + b_6 x_k^4 + b_7 x_k^3 + b_8 x_k^2 + b_9 x_k + b_{10} $

\vspace{1\baselineskip}

where
\[
\begin{pmatrix}
\ b_1 \\
\ b_2 \\
\ b_3 \\
\ b_4\\
\ b_5 \\ 
\ b_6 \\ 
\ b_7 \\ 
\ b_8 \\ 
\ b_9 \\ 
\ b_{10} \\ 
\end{pmatrix}
\approx
\begin{pmatrix}
\ 6.505 \cdot 10^{-12} \\
\ -2.358 \cdot 10^{-10} \\
\-2.009 \cdot 10^{-6} \\
\-8.207 \cdot 10^{-6} \\
\ 6.216 \\ 
\ -2.184 \cdot 10^{-12} \\ 
\ -1.574 \cdot 10^{-08} \\ 
\ -1.790 \cdot 10^{-5} \\ 
\ 3.587 \cdot 10^{-2} \\ 
\ 6.216 \\ 
\end{pmatrix}
\]

\begin{figure}
 \begin{minipage}[b]{.46\linewidth}
  \centering
  \includegraphics[width=\linewidth]{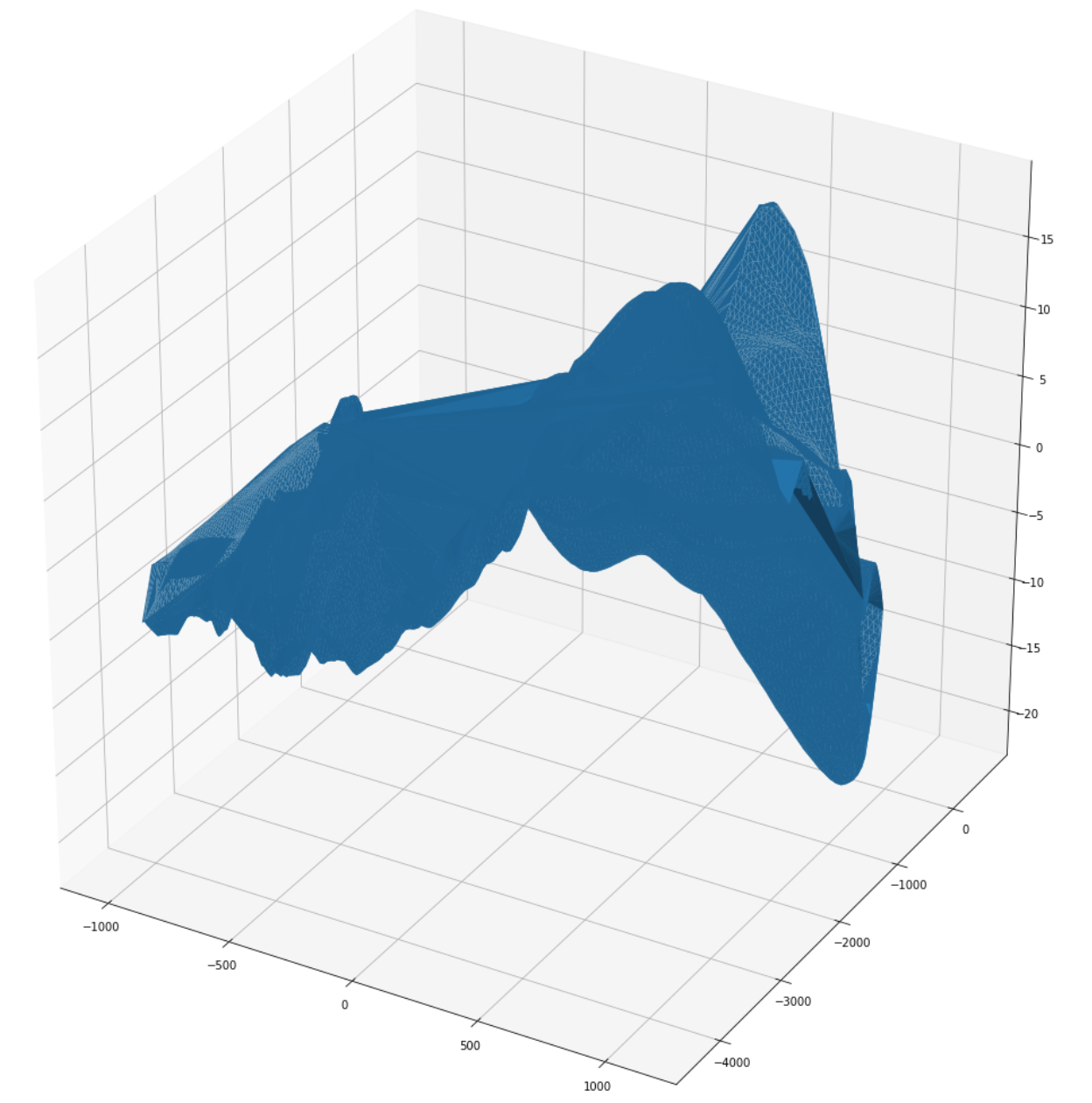}
  \caption{Actual wind force distribution in u direction}
  \label{fig:wind_u_actual}
 \end{minipage} \hfill
 \begin{minipage}[b]{.46\linewidth}
  \centering
  \includegraphics[width=\linewidth]{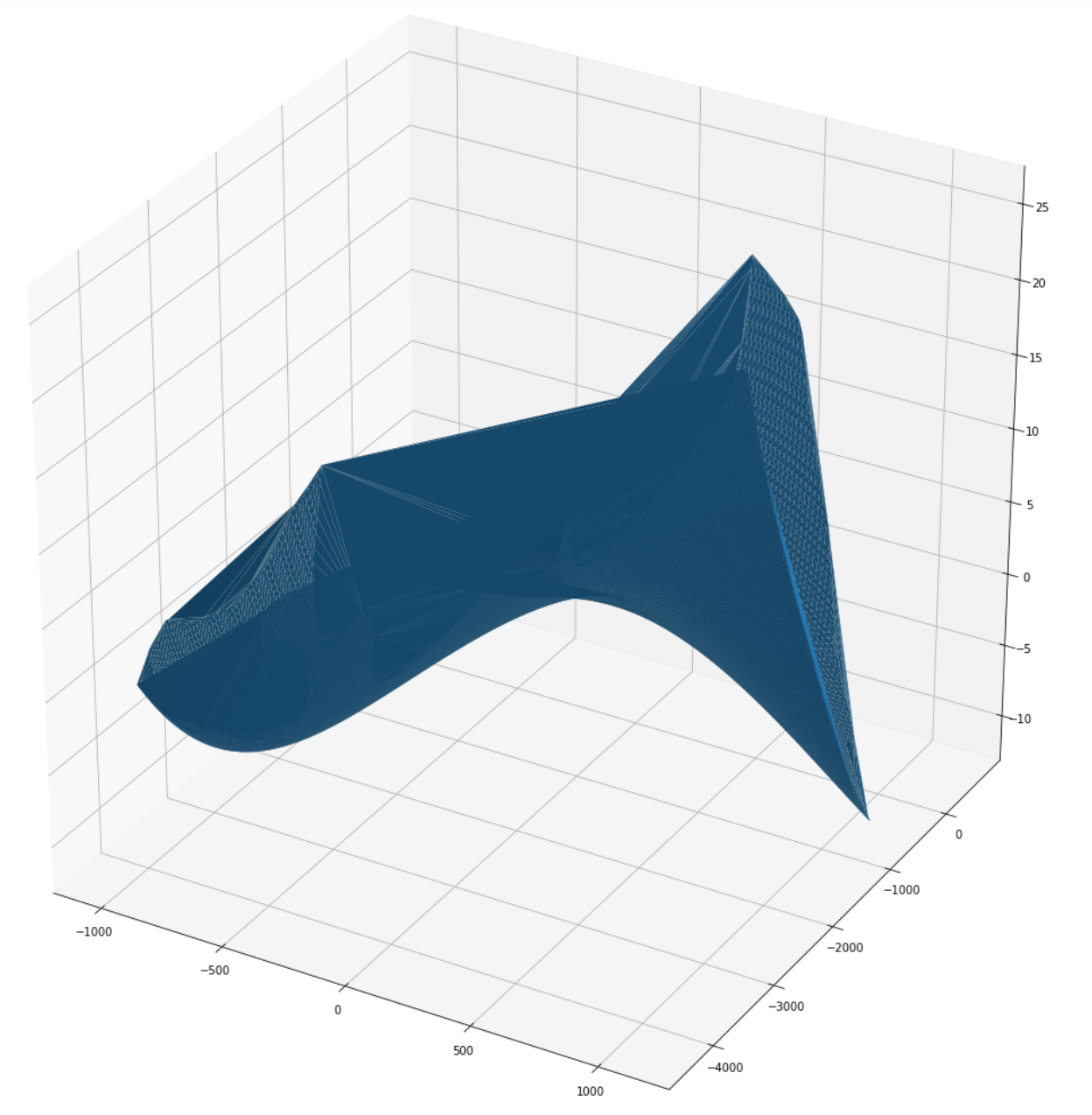}
  \caption{Approximate function in u direction}
  \label{fig:wind_u_approx}
 \end{minipage}
\end{figure}

\begin{figure}
 \begin{minipage}[b]{.46\linewidth}
  \centering
  \includegraphics[width=\linewidth]{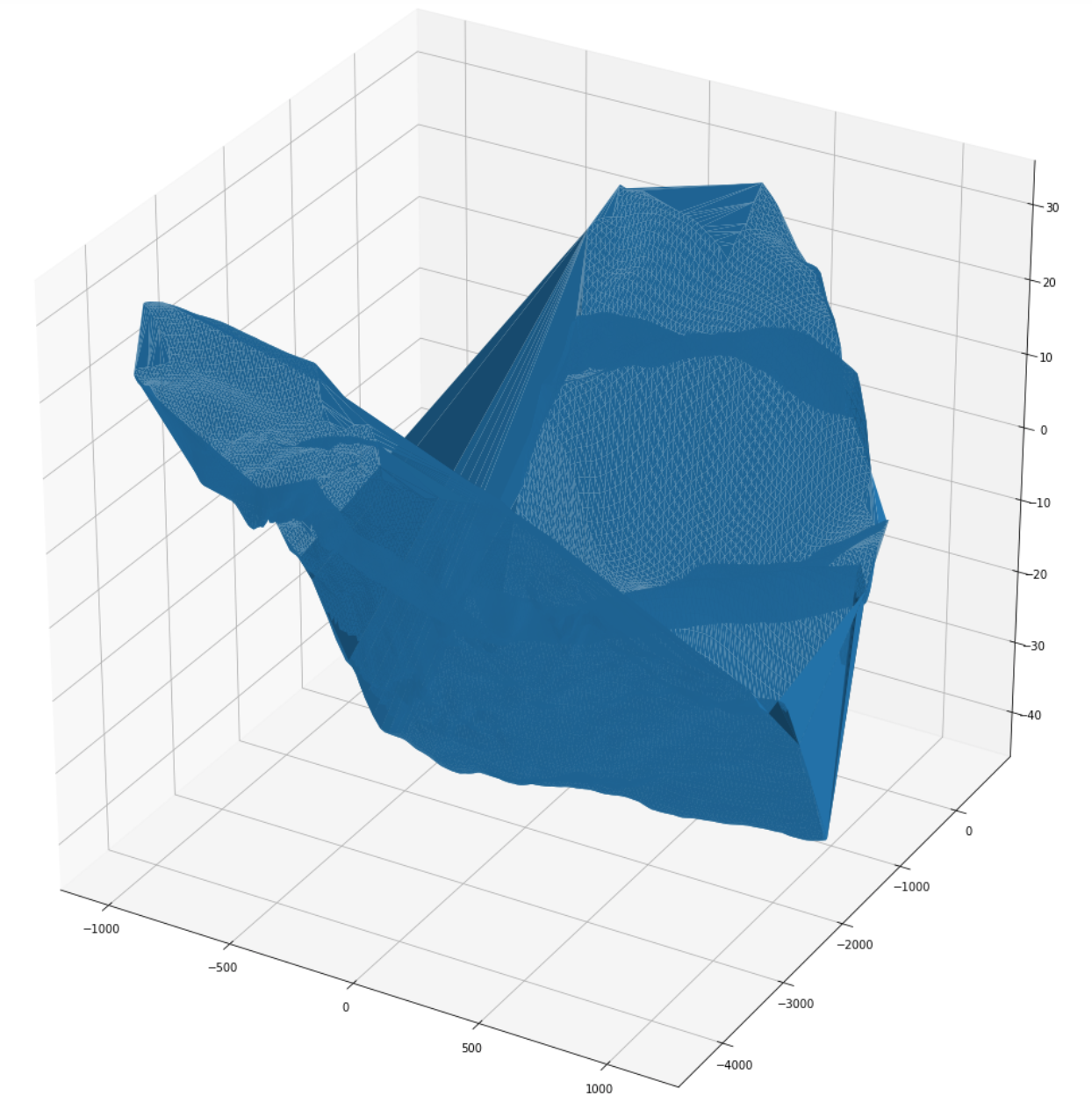}
  \caption{Actual wind force distribution in v direction}
  \label{fig:wind_v_actual}
 \end{minipage} \hfill
 \begin{minipage}[b]{.46\linewidth}
  \centering
  \includegraphics[width=\linewidth]{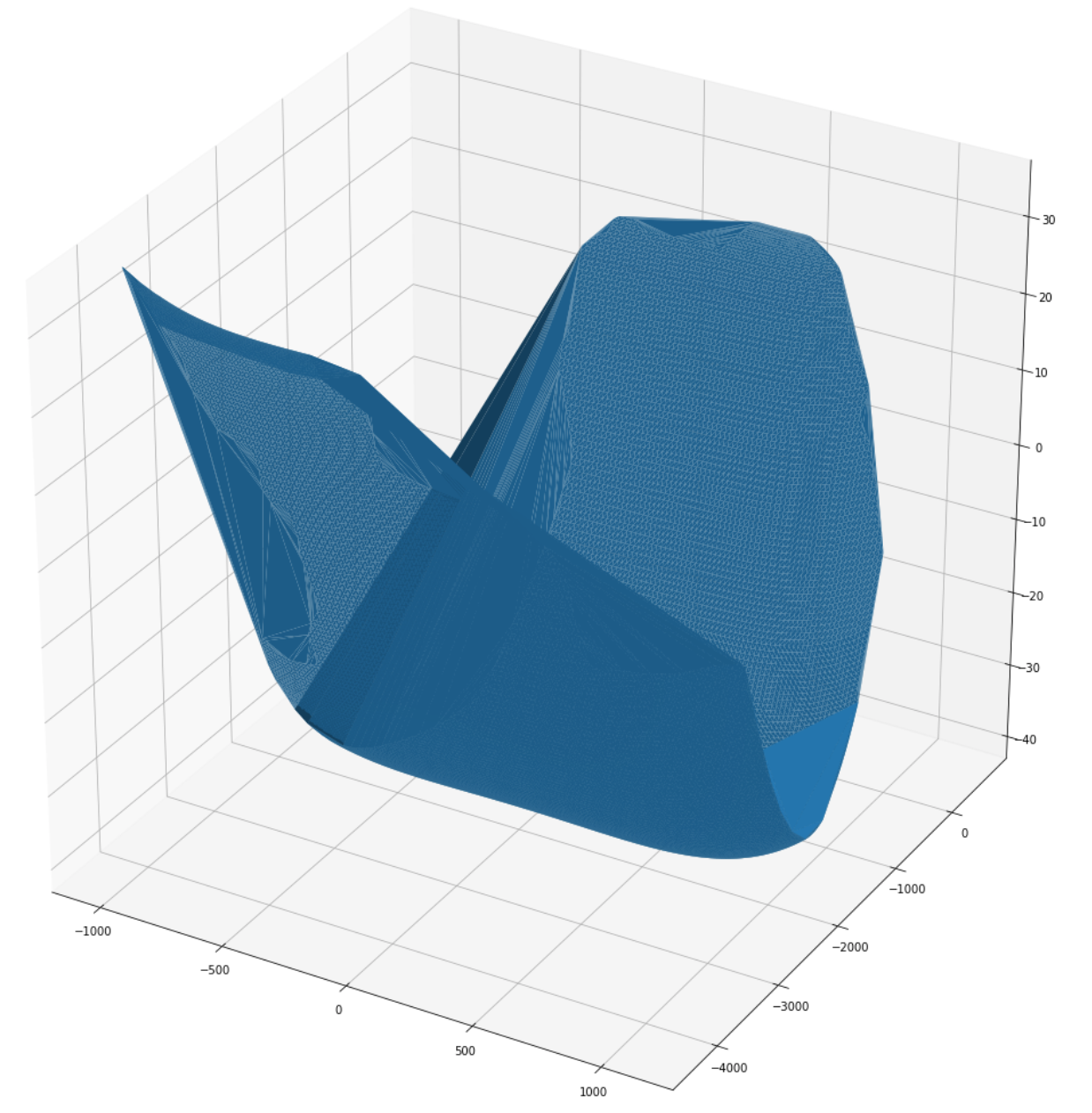}
  \caption{Approximate function in v direction}
  \label{fig:wind_v_approx}
 \end{minipage}
\end{figure}

\subsection{Discretization}
Once the dynamic model, constraints, and the objective function are determined, there are two options to solve it: one is to solve the continuous CFTOC; the other is to solve the discretized CFTOC. The classic way to solve the continuous optimal control problem require to solve the Hamilton–Jacobi–Bellman (HJB) equation, derived from a continuous version of dynamic programming, is a partial differential equation and it is usually hard to solve. It will be even harder when considering the specific constraints.

On the other side, discretized CFTOC is easier to be programmed and there exist solvers that can solve this type
of problem. However, one needs to be careful that discretization can create critical numerical errors, especially in non-linear problems, In this case, we decided to use the discretized CFTOC.

The continuous dynamic model established earlier is therefore discretized. The time derivative of a variable can be approximated by the relationship $x_{k+1} = x_k + \Delta T \cdot \dot{x_k}$. Thus, according to the model above, we finally have the following discretized dynamic relationship: $x_{k+1} = x_k + \Delta T \cdot f(x_k,u_k,w_k)$

\subsection{Non-linear Programming}

\begin{equation*}
\begin{aligned}
& \underset{X,U}{\text{minimize}}
& & \sum_{k=0}^{N-1}  (X_{k} - X_f)^\top Q (X_{k} - X_f) + U_{k}^\top R U_{k} \\
& \text{subject to}
& & X_{lb} \leq X_i \leq X_{ub}, \; i = 0, \ldots, N \\
&&& U_{lb} \leq U_i \leq U_{ub}, \; i = 0, \ldots, N-1 \\
&&& X_{i+1} = X_{i} + \Delta T \cdot f(X_i,U_i,W_i) , \; \forall i \in {N} \\
&&& X_0 = \mathcal{X}_{0} \\
&&& X_N = \mathcal{X}_{f} \\
&&& X_f = \mathcal{X}_{f} \\
\end{aligned}
\end{equation*}

with:
\[
f(X_k,U_k,W_k)
= 
\begin{pmatrix}
\dot{x_k} \\
\dot{y_k} \\
\dot{v_k} \\
\dot{m_k} \\
\dot{\theta_k} \\ 
\end{pmatrix}
=
\begin{pmatrix}
v_k \cdot cos(\theta_k) + w_{xk} \\
v_k \cdot sin(\theta_k) + w_{yk} \\
\frac{2 T_k - C_d \rho A v_{k}^2}{2 m_k} \\
- \eta T_k \\
\varphi_k \\ 
\end{pmatrix}
\]

and

\[
U_k
=
\begin{pmatrix}
T_k \\
\varphi_k
\end{pmatrix}
\]

\section{Results}

The MPC was programmed in MATLAB. The nonlinear solver are very sensitive to the scale of the problem (i.e. numbers of variables). Hence, the total time horizon, time step size, numbers of time step need to be carefully selected.\\

\noindent
$T$ = total time horizon [h]\\
$N$ = total numbers of time step [-]\\
$dT$ = time step size [s]\\

A nonlinear solver IPOPT, which is suitable for large-scale nonlinear optimization, was used. Firstly, we focus on optimizing the flight path without considering wind data and obtain the results as below. To do the optimization, we tried several time by decreasing the whole trip time from 4 hours to 3.6 hours. We wasn't able to make the aircraft arrive at SFO in 3.5 hours. Our optimal travel time is 3.6 hours in total.

\begin{figure}
\begin{center}
    \includegraphics[scale=0.26]{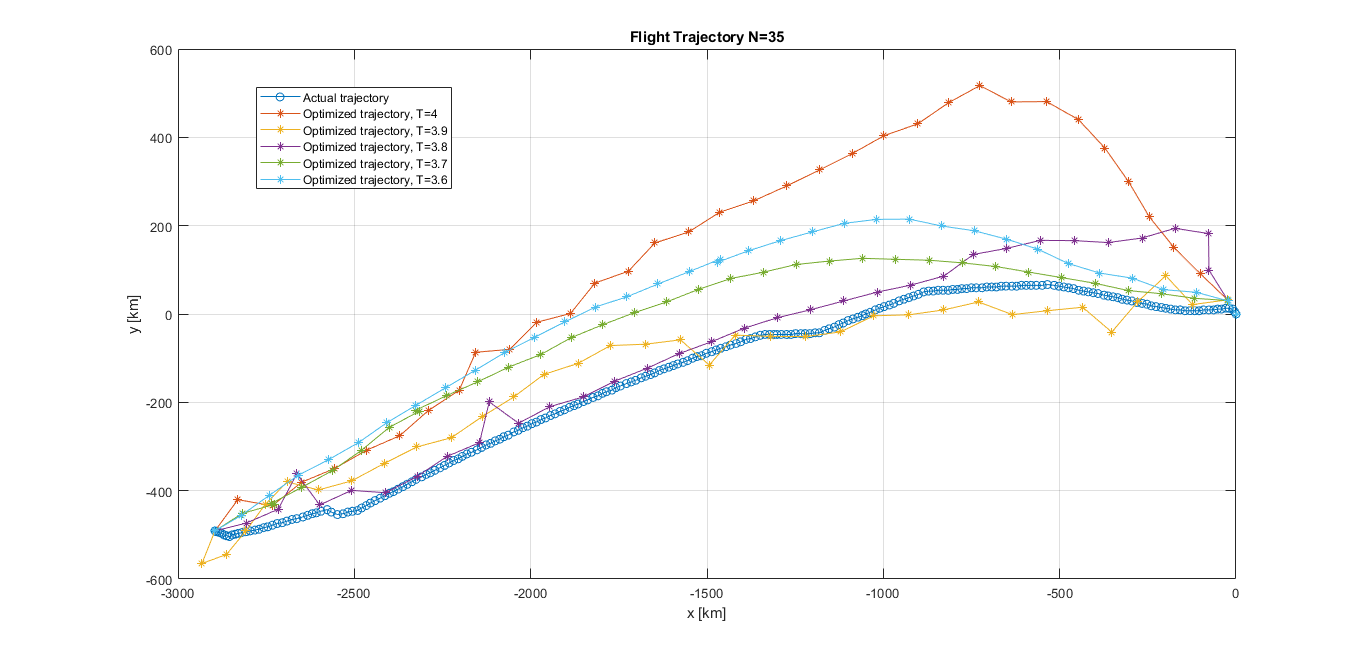}
    \caption{Optimized travel time flight trajectories (no wind) with N=35}
    \end{center}
\end{figure}

Then, based on the polynomial function we fitted for the actual wind data, we apply the function to the state matrix. And the optimal trajectory with wind data is shown as below. If we decrease the total travel time to 3.7 hours, the airplane cannot reach the destination SFO. This time, out optimal flight time is 3.8 hours.

\begin{figure}
\begin{center}
    \includegraphics[scale=0.26]{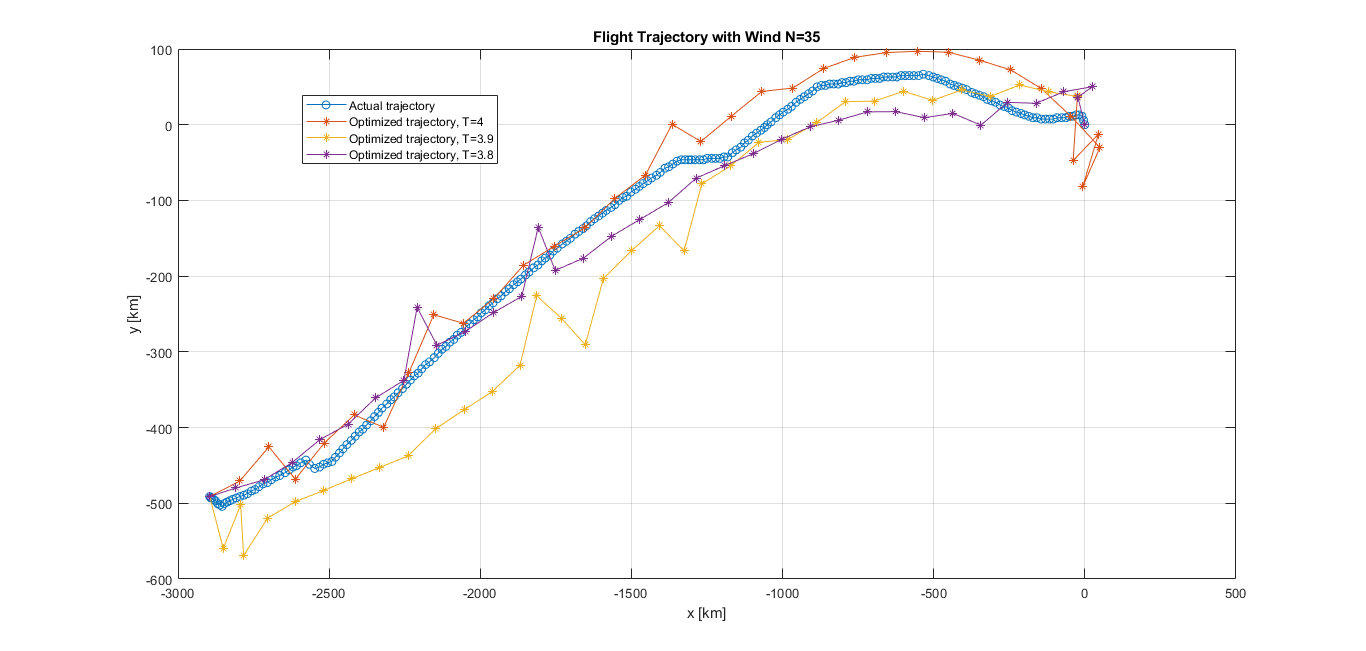}
    \caption{Optimized travel time flight trajectories (with wind) with N=35}
    \end{center}
\end{figure}

Compare the results without wind data and with wind data above, the optimal trajectories with the input of wind data are much closer to the actual flight path, which shows an acceptable result and the fitted wind function performs very well.

\begin{figure}
\begin{center}
    \includegraphics[scale=0.26]{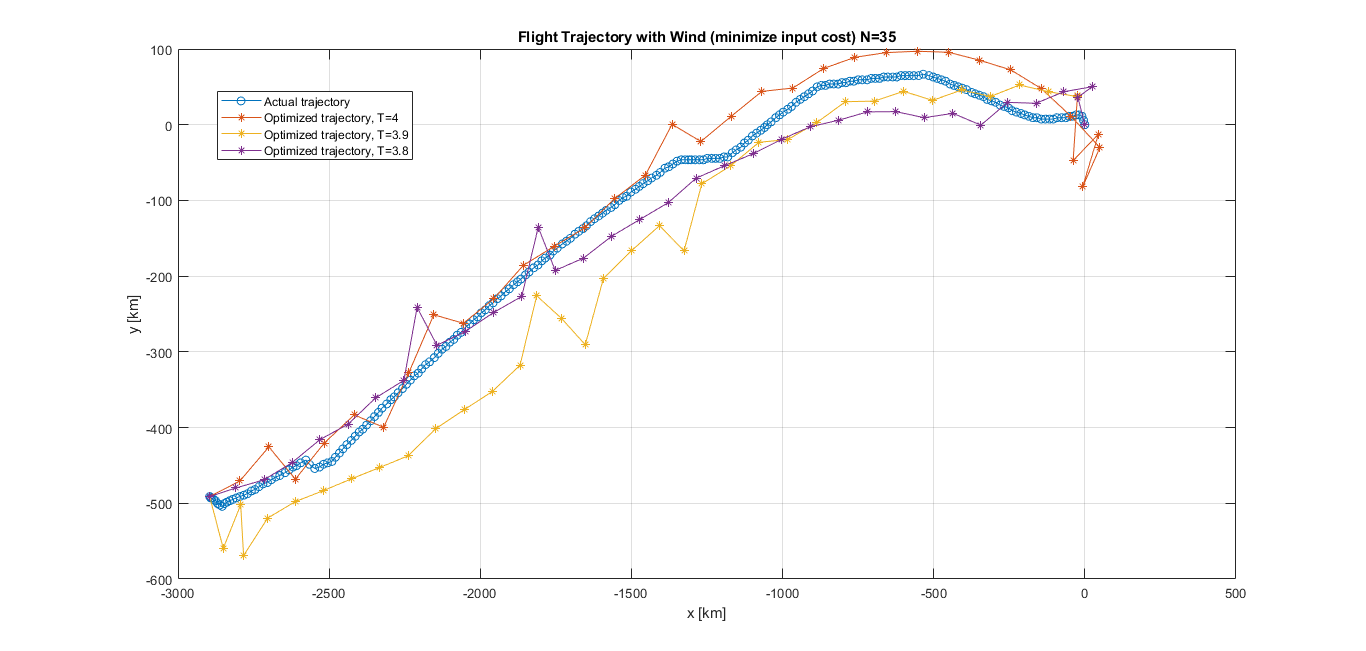}
    \caption{Optimized fuel cost for the flight trajectories (with wind) with N=35}
    \end{center}
\end{figure}

\section{Discussion}
Our project is aimed to minimize the total travel time and fuel cost to obtain the optimal flight path from Chicago O’Hare International Airport (ORD) to San Francisco International Airport (SFO) taking forecast wind data into account. And compare it to the actual flight truck from the existing sample data.
Since the optimization is a nonlinear, non-convex problem, the main challenge is to find a suitable way to solve this optimization problem. The other challenge is that in a nonlinear, non-convex optimization problem, any solution is a local optimum solution, and there is no grantee that it is a global optimum solution.
Although there exist some solvers that can solve nonlinear optimization problem, we found out that these solvers are very sensitive to the constraints (especially to the nonlinear aircraft dynamics), the scale of the problem (i.e. numbers of variables).
For example, for a given total time horizon T=5 [hour], time step =60[s], which lead to total numbers of step = 300. The solver (IPOPT) could not find a solution. However, when reduced the total numbers of step to 40 given T=4 [hour] (i.e. time step = 360[s]), the solver can find a solution.
Also, using different parameters (i.e. different total time horizon, different numbers of step), the solution can change drastically. However, we can see that with certain set up, the optimized path is very similar to the actual path. For a travel time focus optimization, the path is reasonable, since not considering fuel cost, the path associate with the minimum travel time should be a straight line on the surface, which is an arc on a x-y plane.
While taking wind into consideration, the path is even much closer to the actual flight, which indicates that the aircraft model is reasonable and the aircraft performance set up is very close to reality.
The reason might due to nonlinear programming and the defect of the nonlinear solver. Since the behavior of the nonlinear problem is very complicated and hard to handle. And how to handle the nonlinear problem and make a robust optimal control model for this aircraft model should be considered in the future.

\section{Conclusion}
It is possible to use optimal control method to do a flight path optimization with a nonlinear model. However, due to non-linearity of the problem and the problem’s scale, the total number of variables needs to be carefully selected in order to obtain a reasonable result. Hence, the large scale nonlinear problem was reduced to a smaller scale by increase the time step to avoid the solver issue. And this can lead to a low resolution result, which is a drawback of the constrained finite optimal method.
However, comparing the result and the actual flight trajectories and the, the solution is reasonable, which indicate that this optimal control method still have the potential to solve the flight path optimization problem. To obtain a more robust result, further improvement needs to be done in handling the non-linearity of the optimization problem.

\bibliographystyle{IEEEtran}
\bibliography{refs.bib}
\vspace{12pt}

\end{document}